\documentclass{amsart}

\input xy
\xyoption{all}

\usepackage{geometry}
\usepackage{amsmath}
\usepackage{amssymb}
\usepackage{hyperref}
\usepackage{cite}
\usepackage[final]{showkeys}
\usepackage{hyperref}
\usepackage{setspace}
\usepackage{amsthm}  
\newtheorem{same}{This should never appear}[section]
\newtheorem{defin}[same]{Definition}
\newtheorem{claim}[same]{Claim}

\newtheorem{theorem}[same]{Theorem}

\newtheorem{lemma}[same]{Lemma}
\newtheorem{fact}[same]{Fact}

\newtheorem{cor}[same]{Corollary}
\newtheorem{prop}[same]{Proposition}

\newbox\noforkbox \newdimen\forklinewidth
\setbox0\hbox{$\textstyle\smile$}
\setbox1\hbox to \wd0{\hfil\vrule width \forklinewidth depth-2pt
 height 10pt \hfil}
\setbox\noforkbox\hbox{\lower 2pt\box1\lower 2pt\box0\relax}
\def\unionstick{\mathop{\copy\noforkbox}\limits}

\def\nonfork_#1{\unionstick_{\textstyle #1}}

\setbox0\hbox{$\textstyle\smile$}
\setbox1\hbox to \wd0{\hfil{\sl /\/}\hfil}
\setbox2\hbox to \wd0{\hfil\vrule height 10pt depth -2pt width
              \forklinewidth\hfil}
\newbox\doesforkbox
\setbox\doesforkbox\hbox{\lower 2pt\box1 \lower 2pt\box2\lower2pt\box0\relax}
\def\nunionstick{\mathop{\copy\doesforkbox}\limits}

\def\fork_#1{\nunionstick_{\textstyle #1}}

\newcommand{\bL}{\mathbb{L}}

\newcommand{\LS}{\text{LS}}

\newcommand{\ba}{\bold{a}}
\newcommand{\bb}{\bold{b}}
\newcommand{\bc}{\bold{c}}

\newcommand{\bm}{\bold{m}}

\newcommand{\K}{\mathbb{K}}
\newcommand{\bx}{\bold{x}}

\newcommand{\by}{\bold{y}}

\newcommand{\rest}{\upharpoonright}
\newcommand{\id}{\textrm{id}}

\newcommand{\Mod}{\te{Mod }}

\newcommand{\crit}{\text{crit }}
\newcommand{\otp}{\text{otp}}
\newcommand{\df}{\de}

\renewcommand{\P}{\mathcal{P}}
\newcommand{\cP}{\P}

\newcommand{\de}{\text{def}}
\newcommand{\im}{\text{im }}

\newcommand{\te}[1]{\textrm{#1}}

\title{Definable Coherent Ultrapowers and Elementary Extensions}
\author{Will Boney}
\email{wb1011@txstate.edu}
\address{Department of Mathematics, Texas State University, San Marcos, TX, USA}
\date{\today\\This material is based upon work done whilethe author was supported by the National Science Foundation under Grant No. DMS-1402191 and DMS-2339018.}

\begin{document}

\maketitle

\begin{abstract}
We develop the notion of coherent ultrafilters (extenders without normality or well-foundedness).  We then use definable coherent ultraproducts to characterize any extension of a model $M$ in any fragment of $\mathbb{L}_{\infty, \omega}$ that defines Skolem functions by a sufficiently complete (but in $ZFC$) coherent ultrafilter.  We apply this method to various elementary classes and AECs.
\end{abstract}

\tableofcontents

\section{Introduction}

We combine the technique of definable ultrapowers and coherent ultrafilters (a $ZFC$ version of extenders) to characterize extensions of a model in some fragment of $\mathbb{L}_{\infty, \omega}$ (including first order) with definable Skolem functions.  Coherent ultrafilters are a weakening of extenders, a tool from set theory
.  The goal of extenders is to capture some elementary embedding $j:V \to \mathcal{M}$ where $\mathcal{M}$ is an inner model that has some rank-bounded similarity to $V$ by creating a coherent system of ultrafilters that satisfy additional properties ($\kappa$-completeness, normality, and well-foundedness; see Kanamori \cite[Chapter 26]{kanamori} for a reference).  The use of a coherent system of ultrafilters is necessary to separate out the property of `closure under $\lambda$-sequences' (a saturation-like property) from the property of containing, e. g., some $V_\lambda$.  Rather than taking an ultrapower of $V$, the coherence gives a directed system of ultrapowers of $V$ and the ``extender power" (or coherent ultrapower in our terms) is the colimit\footnote{Or `directed limit'} of this system.

Definable ultrapowers are a technique that goes back to Skolem \cite{s-skolemult} (although we use \cite[Chapter 32]{k-lw1w} as our more accessible reference).  Rather than including \emph{all} functions from the index set to the model $M$, the definable ultrapower uses some definable subset of $M$ as the index and only considers the definable functions from that definable subset to the model.  This has the drawback of requiring more structure on $M$, in particular, requiring $Th(M)$ to have definable Skolem functions.  The pay-off, however, is that there is much greater control over what the definable ultrapower looks like.  The work in \cite[Chapter 32]{k-lw1w} exploits this to, from $M \prec N$, create a $ZFC$ ultrafilter that is complete enough to ensure not just first-order elementarity of the extension by the definable ultrapower, but elementarity according to whatever fragment of $\mathbb{L}_{\infty, \omega}$ that captures the relation between $M$ and $N$.  The key result that we generalize is the following (\cite{k-lw1w} deals exclusively with $\bL_{\omega_1, \omega}$, but the adaptation is immediate).

\begin{fact}[{\cite[Chapter 32]{k-lw1w}}]\label{motiv-fact}
Let $\mathcal{F}$ be a fragment of $\mathbb{L}_{\infty, \omega}(\tau)$ and $\psi \in \mathcal{F}$ have definable Skolem functions.  Given $M \prec_\mathcal{F} N$ that satisfy $\psi$, there is an $\mathcal{F}$-complete ultrafilter $U$ on $M$ such the definable ultrapower $\prod^{\df} M/U$ $\mathcal{F}$-elementary embeds into $N$ over $M$.  Moreover, given some $c \in N$, $U$ can be constructed so that $c$ is in the image of this embedding.
\end{fact}

The moreover clause of Fact \ref{motiv-fact} is the key to our contribution.  The coherent ultrafilters allow us to move from specifying some $c \in N$ in the target model to specifying an arbitary subset of the target model.  This allows us to make $\prod^{\df} M/U$ isomorphic to $N$ (see Theorem \ref{main-thm}).  In order to do this, we must show that the $\mathcal{F}$-complete coherent filter derived from $M \prec_\mathcal{F} N$ can be extended to a $\mathcal{F}$-complete coherent ultrafilter; this is the goal of Section \ref{cohult-sec}. In Section \ref{ex-sec}, we work through some examples of classes and apply this theorem.

Note that, in the context of $(V, \in)$ every function $f:\kappa \to V$ is definable with parameters\footnote{By the formula ``$x = f$."}, so the distinction between $\prod V/E$ and $\prod^{\df} V/E$ disappears (when the index is a set).  

After this paper was first circulated, Enayat, Kauffman, and McKenzie circulated \cite{ekm-iterated}, which also deals with coherent ultraproducts (although they call them `dimensional ultrapowers').  Their work focuses on using iterated ultrapowers in the first-order context to build tight indiscernibles without explicit use of the Compactness Theorem or Ramsey's Theorem, which provides insight into automorphisms of models.  \cite[Section 4]{ekm-iterated} also contains a nice overview of the history and applications for the interested reader.  We also thank the anonymous referee for an extremely helpful report with many suggestions and corrections (especially making us aware of Keisler's problem and its relation to Section \ref{ex-sec}.

\section{Background}
A few notes on notation are in order.  We refer to a system of ultrafilters that is coherent (see Definition \ref{main-def}) as a ``coherent ultrafilter" rather than the clunkier but more accurate ``coherent system of ultrafilters."  We also discuss some general facts about coherent ultraproducts.

The motivation and methods of this paper blend together ideas from elementary model theory, infinitary model theory, and large cardinal set theory.  In order to make things more clear, we give a background on some of these ideas.  The results appear as folklore and/or in textbooks, so we forego specific citations; the material can be found in places like \cite[Chapter 6]{changkeisler} and \cite[Chapters 17 and 20]{jech}.

Working in first-order logic, one entry point into the world of systems of ultrapowers is the problem of complete extensions.  We write $M^\#$ for the full expansion of $M$, that is, the expansion to a language that has a function symbol for every function $M^n \to M$ and a relation symbol for every element of $\cP(M^n)$. Given structures $M$ and $N$, we say that $N$ is a \emph{complete extension of $N$} iff there is an expansion $N^+$ of $N$ to the full language of $M$ such that
$$M^\# \prec N^+$$
Any ultrapower of $M$ is a complete extension since ultraproducts commute with restriction of the language.  Moreover, any increasing union of ultrapowers or more complicated system will share this property; these systems are sometimes called \emph{limit ultrapowers}, intimately related to our coherent systems of ultrapowers.  In fact, these limit ultrapowers provide a complete characterization of complete extensions.

\begin{fact}
    Fix two structures $M \prec N$.  $N$ is a complete extension of $M$ iff $N$ is isomorphic to a limit ultrapower of $N$.
\end{fact}

Another appearance of these systems of ultrapowers is in the context of definable Skolem functions. Fix a first-order theory $T$ with definable Skolem functions and models $M \prec N$.  Given any $\ba \in N$, we will write $M(\ba)$ for the closure of $M\cup\{\ba\}$ under the definable functions of $N$; because $T$ has definable Skolem functions, this gives
$$M \prec M(\ba) \prec N$$

Moreover, given two tuples $\ba \subset \bb \in N$, we have $M(\ba) \prec M(\bb)$.  Thus we can write $N$ as the directed colimit of the system 
$$\left\{M(\ba) : \ba \in N\right\}$$
These models $M(\ba)$ can be viewed as definable ultrapowers.  Given $\ba \in N$, we can define the filter of definable subsets of $M$ that $\ba$ satisfies, that is,
\begin{eqnarray*}
    \phi(M)&:=&\{\bm \in M^{\ell(\bx)}: M\vDash \phi(\bm)\}  \text{ for a formula } \phi(\bx)\\
    U_\ba&:=& \left\{\phi(M): N\vDash\phi(\ba)\right\}
\end{eqnarray*}
Then we set $\prod^{\df}_n M$ to be the collection of $M$-definable functions from some $M^n\to M$.  When we mod out by the filters $U_\ba$, we get the definable Skolem ultrapower
$$\prod^\df_nM/U_\ba$$

These Skolem ultrapowers reproduce the models $M(\ba)$ defined above:
\begin{fact}
Let $T$ be a theory with definable Skolem functions and $M\prec N$ model $T$.
    \begin{enumerate}
        \item Given $\ba\in N^n$, we have
        $$M(\ba) \cong \prod^\df_nM/U_\ba$$
        \item $N$ is the directed colimit of Skolem ultrapowers of $M$.
    \end{enumerate}
\end{fact}
Keisler \cite[Chapter 32]{k-lw1w} explored the definable ultrapowers in the context of fragments of $\bL_{\omega_1, \omega}$ for $n=1$.  We extend this analysis to the full result.

These results have analogues in set theory to full ultrapowers to stronger logics.  Having a countably complete ultrapower already implies the existence of a measurable cardinals, so large cardinals are immediately involved.

For single ultrafilters, given a non-identity elementary embedding\footnote{This is the common way of writing $V \precneqq M$ in the set-theoretic framework} $j:V \to M$, set $\kappa = \crit j = \min\{\alpha \in \text{Ord}:\alpha \neq j(\alpha)\}$.  Then one can define an filter as above
$$U_\kappa = \{X \subset \kappa:\kappa \in j(X)\}$$
This is in fact a $\kappa$-complete, normal ultrafilter, giving rise to an $\bL_{\kappa, \kappa}$-elementary embedding
$$j_U:V \to \prod V/U_\kappa$$
Further, this ultrapower acts like a version of $V(\kappa)$, in the sense that there is an embedding $k:\prod V/U_\kappa \to M$ such that $\kappa \in \im k$ and $j = k \circ j_U$.

Stronger large cardinals (especially strong) are characterized by building a stronger embeddings that come from systems of ultrafilters.  These systems are called \emph{extenders}.  Namely, a strong embedding is one of the form $j:V\to M$ with $\crit j = \kappa$ with $V_\lambda \subset M$ for $\lambda >\kappa$.  For each $\ba \in V_\lambda$, one can define a $\kappa$-complete ultrafilter $U_\ba$ on $\kappa^{\ell(\ba)}$ as before.  Then, similar to the definable ultrapower case, the ultrapowers
$$\left\{\prod V/U_{\ba} : \ba \in  V_\lambda\right\}$$
form a directed system and we can form the directed colimit
$$j_0:V \to M_0$$
We cannot expect to capture all of $M$ in this way because $M$ is a class and we must work with a set of ultrafilters.  Nonetheless, we can find an embedding $k:M_0 \to M$ with the properties that $V_\lambda \subset \im k$ and $j = k\circ j_0$.

Our coherent systems of ultrafilters are based on these extenders after stripping away the extra large cardinal properties.

\section{Coherent Ultrafilters} \label{cohult-sec}

The following notion of coherent ultrafilter is designed to take the key model-theoretic features of an extender (see \cite[Chapter 26]{kanamori}) that allow the construction of $\prod V/E$, and remove the parts that correspond to beyond-ZFC strength.  We allow the arity $\mu$ of the extender to be larger than $\omega$ (that is, a $(\kappa, \lambda)$-extender is a $(\omega, \kappa, \lambda)$-coherent ultrafilter with extra structure), although we will not use it here.  Also, we remove the requirement that $\kappa$ be a cardinal, replacing it with an arbitrary set $A$.  Although this does not matter for the full coherent ultrapower since there is a bijection between $A$ and $|A|$, it is an important consideration for the definable coherent ultrapower since the bijection is rarely definable.  

We use a slightly different formalism for coherent filters than typical.  Normally, for $\ba \in [\lambda]^{<\omega}$, one of the following is used:
\begin{itemize}
	\item $F_\ba$ is a filter on $[A]^{|\ba|}$, so there is a canonical pairing between $\ba$ and any $s \in [A]^{|\ba|}$ that pairs each element of $\ba$ with the corresponding one of $s$ (this is used in \cite{kanamori})
	\item $F_\ba$ is a filter on ${}^\ba A$, so, for $s \in {}^\ba A$, the pairing of elements is given by the function (this is used in \cite{ms-extenders})
\end{itemize}
The pairing is necessary to define the projections $\pi^\bb_\ba$ for $\ba \subset \bb$, which is how one makes sense of coherence.  Since we want to work with definable objects in some model, neither of these are available.  Instead, we will take $\ba \in [\lambda]^{<\omega}$ and set $F_\ba$ to be a filter on ${}^{|\ba|} A$.  Then we can pair the $i$th element of $\ba$ (under the ordering inherited from the ordinals) with the $i$th element of $s \in {}^{|\ba|} A$ in the order inherited from the function.  This is clunkier, but is necessitated by our situation.  As part of this, $s \in {}^n A$ will be considered both as an $n$-tuple of elements and as a function from $n$ to $A$.  When we move to longer $\ba$'s, we will use ${}^{\otp(\ba)}A$ to make use of the ordering.

\begin{defin} \label{main-def}
Let $A$ be a set, $\mu\leq \lambda$ cardinals.  
\begin{enumerate}
    \item $[A]^{\mu}$ is the collection of $\mu$-sized subsets of $A$ and $(\lambda)^{\mu}$ are the subsets of $\lambda$ with order-type $\mu$.  Then
    \begin{eqnarray*}
        [A]^{<\mu} &:=& \bigcup_{\alpha < \mu} [A]^\alpha\\
        (A)^{<\mu} &:=& \bigcup_{\alpha < \mu} (A)^\alpha
    \end{eqnarray*}
	\item Let $\alpha < \beta < \mu$, $\ba \in (\lambda)^{\alpha}$, and $\bb \in (\lambda)^\beta$ with $\ba \subset \bb$.  Set $p_{\ba, \bb}:\otp(\ba) \to \otp(\bb)$ be the unique order-preserving injection such that the $i$th member of $\ba$ is the $p_{\ba, \bb}(i)$th member of $\bb$.  Set $\pi^\bb_{\ba}: {}^\beta A \to {}^\alpha A$ by 
	$$\pi^\bb_{\ba}(s) = s \circ p_{\ba, \bb}$$
    \item Given $\ba\subset \bb$ and $s \in \otp(\ba)$, set
    $$s\uparrow\bb:=s \circ p_{\ba,\bb}^{-1}$$
    \item If $\ba \subset \bb \subset \bc$ and $X \subset {}^{\otp(\bb)} A$, then we use the notation
    \begin{eqnarray*}
        \pi^{\bb}_{\ba}X&:=&\pi^{\bb}_{\ba}[X]= \{\pi^{\bb}_{\ba}(s): s \in X\}\subset {}^{\otp(\ba)}A\\
        (\pi^{\bc}_{\bb})^{-1} X &=&(\pi^{\bc}_{\bb})^{-1}[X]= \left\{ s : \pi^\bc_\bb(s) \in X\right\} \subset {}^{\otp(\bc)} A
    \end{eqnarray*}
	\item We say that $F$ is a \emph{$(\mu, A, \lambda)$-coherent filter} iff $F = \{F_\ba \mid \ba \in (\lambda)^{<\mu}\}$ such that
	\begin{enumerate}
		\item each $F_\ba$ is a filter on ${}^{\otp(\ba)} A$; and
		\item if $\ba \subset \bb$, then for any $X \subset {}^{\otp(\ba)} A$, we have
		$$X \in F_\ba \iff \{s \in {}^{\otp(\bb)} A \mid \pi^\bb_\ba(s) \in X\} \in F_\bb$$
	\end{enumerate}
	\item We say that $E$ is a \emph{$(\mu, A, \lambda)$-coherent ultrafilter} iff it is a $(\mu, A, \lambda)$-coherent filter with each $E_\ba$ an ultrafilter.
	\item A $(\mu, A, \lambda)$-coherent filter $F$ is \emph{proper} iff no $F_\ba$ contains the empty set.
	\item Given two $(\mu, A, \lambda)$-coherent filters $F$ and $F^*$, we say that \emph{$F^*$ extends $F$} or $F \Subset F^*$ iff $F_\ba \subseteq F^*_\ba$ for all $\ba \in (\lambda)^{<\mu}$.
	\item Given $\bb \subset \bc$ and $X \subset {}^{\otp(\bc)} A$, we say that $X$ is full over $\bb$ iff for every $s, t \in {}^{\otp(\bc)} A$ such that $\pi^\bc_\bb(s) = \pi^\bc_\bb(t)$, we have
	$$s \in X \iff t \in X$$
	\item If we omit $\mu$ from any of the above definitions, we mean $\mu = \omega$.

\end{enumerate}
\end{defin}
We deal almost exclusively with the case $\mu=\omega$.  In this case, we write $[\lambda]^{<\omega}$ instead of $(\lambda)^{<\omega}$ and $|\ba|$ instead of $\otp(\ba)$; they are the same, and the former follows convention.  We could also vary this definition to allow different $F_{\{a\}}$ to live on different sets, rather than the same one; this is sometimes necessary with extenders.  However, since we don't impose extra conditions like normality, we can extend the underlying set so all $F_{\{a\}}$ live on the same one in the context of coherent filters.

  Fullness is a useful notion because it allows us to write the coherence condition as the conjunction of the following:
\begin{itemize}
	\item If $\ba\subset \bb$ and $X \in F_\ba$, then $(\pi^\bb_\ba)^{-1}X \in F_\bb$.
	\item If $\ba \subset \bb$ and $X \in F_\bb$ is full over $\ba$, then $\pi^\bb_\ba X \in F_\ba$.
\end{itemize}

This is used in, e.g., Lemma \ref{onestep-lem}.
There are two useful constructions that we will use throughout this section, both about extending certain objects.

For the first, we show that we can extend sets to minimal full extensions.

\begin{lemma}
Given $\ba \subset \bb$, each subset $X$ of ${}^{\otp(\bb)} A$ has a minimal extension that is full over $\ba$; we denote this set $X^+$ and call it the fulfillment of $X$ over $\ba$.
\begin{enumerate}
    \item This set has the property that $\pi^\bb_\ba X = \pi^\bb_\ba X^+$.
    \item If $X$ is a filter, then so is $X^+$.
\end{enumerate}
\end{lemma}

{\bf Proof:} Set
$$X^+:=\left(\pi^\bb_\ba\right)^{-1}\pi^\bb_\ba X = \{s \in {}^{\otp(\bb)}A: \exists t \in X. \pi^{\bb}_{\ba}(s) = \pi^{\bb}_{\ba}(t)\}$$
Then straightforward to verify.\hfill \dag\\

We will also often be in a situation where we have $\ba \in (\lambda)^{<\mu}$ and $s, t, \dots$ that are defined for subsets of $\ba$, and we want to find a mutual extension of them to ${}^{\otp(\ba)}A$.  Using our notation, $\pi^{\bb}_{\ba}(s)$ serves as the restriction of $s \in {}^{\otp(\bb)}A$ to $\otp(\ba)$.  Dually, $s \uparrow \bc$ is the copy of $s\in {}^{\otp(\bb)}A$ with domain the proper subset of $\otp(\bc)$.

Then we can formulate the necessary and sufficient condition for finding extensions.

\begin{lemma}\label{ext-lem}
    Let $\ba, \bb, \bc \in (\lambda)^{<\mu}$ and $s \in {}^{\otp(\ba)}A$, $t\in {}^{\otp(\bb)}A$ with $\ba \subset \bc$ and $\bb \subset \bc$.  Then the following are equivalent:
    \begin{enumerate}
        \item there is an extension in ${}^{\otp(\bc)}A$ of $s \uparrow \bc$ and $t \uparrow \bc$.
        \item $$\pi^{\ba}_{\ba \cap \bb}(s) = \pi^{\bb}_{\ba \cap \bb}(t)$$
    \end{enumerate}
\end{lemma}

We won't use it, but this has a natural extension to any number of functions we want to extend.

The goal of this section is to generalize the fact that any nonprincipal filter can be extended to a nonprincipal ultrafilter to the context of coherent (systems of) filters and ultrafilters; this is Theorem \ref{cohultext-thm}.  This seems trickier than expected since adding a new set to some $F_\bb$ commits new sets to the rest of the filters (by coherence) and could break the nontriviality of some $F_\ba$.  However, as the following lemmas show, the original coherence is enough to ensure this doesn't happen.

The first lemma is a technical calculation that will come in handy later.

\begin{lemma}\label{nicefull-lem}
Let $\ba, \bb, \bc \in [\lambda]^{<\mu}$ with $\bb \subset \bc$, $X \subset {}^{\otp(\bc)} A$, and $Y \subset {}^{\otp(\ba)} A$ such that $X$ is full over $\bb$ and 
$$(\pi^\bc_{\ba\cap \bc})^{-1}\pi^\ba_{\ba\cap \bc} Y \subset X$$
Set $Y^+$ to be the fulfillment of $Y$ over $\ba \cap \bb$.  Then
\begin{enumerate}
	\item $(\pi^\bc_{\ba\cap \bc})^{-1} \pi^\ba_{\ba\cap \bc} Y^+ \subset X$; and
	\item $\pi^\ba_{\ba\cap \bc} Y^+$ is full over $\ba \cap \bb$.
\end{enumerate}
\end{lemma}

{\bf Proof:} For (1), let $s \in (\pi^\bc_{\ba\cap \bc})^{-1} \pi^\ba_{\ba\cap \bc} Y^+$.  Unraveling the definitions (and using $\bb \subset \bc$, this means that there is $t \in Y$ such that
$$\pi^\bc_{\ba \cap \bb}(s) = \pi^{\ba}_{\ba \cap \bb}(t)$$
Then, we also have
$$\pi^{\bb}_{\ba \cap \bb} \left(\pi^{\bc}_{\bb}(s)\right) = \pi^{\bc}_{\ba\cap\bb}(s) = \pi^{\ba}_{\ba\cap \bb}(t) = \pi^{\ba\cap\bc}_{\ba\cap\bb}\left(\pi^{\ba}_{\ba\cap\bc}(t)\right)$$
By Lemma \ref{ext-lem}, this guarantees tha there is $s' \in {}^{\otp(\bc)}A$ that extends 
$$\left(\pi^{\bc}_{\bb}(s)\right)\uparrow \bc \text{ and }\left(\pi^{\ba}_{\ba \cap \bc}(t) \right) \uparrow \bc$$
Then, 
$$\pi^{\bc}_{\ba\cap\bc}(s') =\pi^{\ba}_{\ba\cap\bc}(t) \in \pi^{\ba}_{\ba\cap\bc}Y$$
Thus $s' \in X$ by the assumption.  Also, $\pi^{\bc}_{\bb}(s')=\pi^{\bc}_{\bb}(s)$, so we also have $s \in X$ as desired by the fullness over $\bb$.

For (2), to show the fullness, let $s \in \pi^{\ba}_{\ba\cap\bc}Y^+$ and $s' \in {}^{\otp(\ba\cap\bc)} A$ such that
$$\pi^{\ba\cap\bc}_{\ba\cap\bb}(s) = \pi^{\ba\cap\bc}_{\ba\cap\bb}(s')$$
Unraveling the definition of $s$, there is $t \in Y^+$ such that $s =\pi^{\ba}_{\ba\cap\bc}(t)$.  Find $t' \in {}^{\otp(\ba)}A$ extending $s'\uparrow\ba$.  Then
$$\pi^{\ba}_{\ba\cap\bb}(t') = \pi^{\ba\cap\bc}_{\ba\cap\bb}(s') = \pi^{\ba}_{\ba\cap\bb}(t)$$
Since $Y^+$ is full over $\ba\cap\bb$, we have $t'\in Y^+$.  Then
$$s' =\pi^{\ba}_{\ba\cap\bc}(t') = \pi^{\ba}_{\ba\cap\bc}Y^+$$
\hfill\dag\\

The following lemma tells us that it's enough to push new sets down and then up.

\begin{lemma}\label{du=ud-lem}
Suppose $\ba, \bb, \bc \in (\lambda)^{<\mu}$ such that $\ba \subset \bc$ and $\bb \subset \bc$.  Let $X \subset {}^{\otp(\ba)} A$ be full over $\ba \cap \bb$.  Then
$$(\pi^\bb_{\ba\cap \bb})^{-1}\pi^\ba_{\ba \cap \bb} X = \pi^\bc_\bb (\pi^\bc_\ba)^{-1} X$$
\end{lemma}

{\bf Proof:} In one direction, let $s \in (\pi^{\bb}_{\ba\cap\bb})^{-1}\pi^{\ba}_{\ba\cap\bb}X$.  This means there is $t \in X$ such that $\pi^{\bb}_{\ba\cap\bb}(s) = \pi^{\ba}_{\ba\cap\bb}(t)$.  By Lemma \ref{ext-lem}, there is $s'\in {}^{\otp(\bc)}A$ extending $s\uparrow\bc$ and $t\uparrow \bc$.  Then 
$$\pi^\bc_\ba(s')=t \in X$$
so $s' \in (\pi^{\bc}_{\ba})^{-1}X$.  Also, $s = \pi^{\bc}_{\bb}(s')$, so $s \in \pi^{\bc}_{\bb}(\pi^{\bc}_{\ba})^{-1}X$, as desired.

In the other direction, let $s \in \pi^{\bc}_{\bb}(\pi^{\bc}_{\ba})^{-1}X$.  This means that there is $t \in {}^{\otp(\bc)}A$ such that
$$\pi^{\bc}_{\bb}(t)=s\text{ and }\pi^{\bc}_{\ba}(t) \in X$$
Then
$$\pi^{\bb}_{\ba\cap\bb}(s) = \pi^{\bc}_{\ba\cap\bb}(t) = \pi^{\ba}_{\ba\cap\bb}\left(\pi^{\bc}_{\ba}(t)\right)$$
So $s \in (\pi^{\bb}_{\ba\cap\bb})^{-1}\pi^{\ba}_{\ba\cap\bb}X$
\hfill\dag\\

\begin{lemma}\label{extfip-lem}
Suppose $F$ is a $(\mu, A, \lambda)$-coherent filter and $F^*_\ba$ is a proper filter on ${}^{\otp(\ba)} A$ extending $F_\ba$.  Then, for all $\bb \in (\lambda)^{<\mu}$, 
$$F_\bb \cup \{ (\pi^\bb_{\ba \cap \bb})^{-1} \pi^\ba_{\ba \cap \bb} X \mid X \in F^*_\ba\text{ is full over }\ba \cap \bb\}$$
has the finite intersection property.
\end{lemma}

{\bf Proof:}  Suppose that this is false for some $\bb \in (\lambda)^{<\mu}$.  Then there are $Y \in F_\bb$ and $X_1, \dots, X_n \in F^*_\ba$ that are full over $\ba \cap \bb$ such that
$$Y \cap \bigcap_i (\pi^\bb_{\ba \cap \bb})^{-1} \pi^\ba_{\ba \cap \bb} X_i = \emptyset$$
Set $X = \cap_i X_i \in F^*_\ba$.  Then we have
$$Y \cap (\pi^{\bb}_{\ba \cap \bb})^{-1} \pi^\ba_{\ba \cap \bb} X = \emptyset$$
Set $Y^+$ to be the fulfillment of $Y$ over $\ba \cap \bb$.  Then $Y^+ \in F_\bb$ and $(\pi^\ba_{\ba \cap \bb})^{-1} \pi^\bb_{\ba \cap \bb} Y^+ \in F_\ba$.  Since $F^*_\ba$ is proper, there is 
$$t \in (\pi^{\ba}_{\ba \cap \bb})^{-1} \pi^\bb_{\ba \cap \bb} Y^+ \cap X$$
Thus $t \in X$ and $\pi^\ba_{\ba \cap \bb}(t) = \pi^\bb_{\ba \cap \bb}(t_0)$ for some $t_0 \in Y^+$.  By definition of $Y^+$, there is $t_1 \in Y$ such that $\pi^\ba_{\ba \cap \bb}(t) = \pi^\bb_{\ba \cap \bb}(t_1)$.  This means that $t_1 \in (\pi^\bb_{\ba \cap \bb})^{-1}\pi^\ba_{\ba \cap \bb} X$.

Thus, $t_1 \in Y \cap (\pi^\bb_{\ba \cap \bb})^{-1} \pi^\ba_{\ba \cap \bb} X$, a contradiction.\hfill \dag\\

\begin{lemma}\label{onestep-lem}
Let $F$ be a $(\mu, A, \lambda)$-coherent filter and $F^*_\ba$ a proper filter on ${}^{\otp(\ba)} A$ extending $F_\ba$.  For each $\bb \in (\lambda)^{<\mu}$, set $F^*_\bb$ to be the filter generated by 
$$F_\bb \cup \{ (\pi^\bb_{\ba \cap \bb})^{-1} \pi^\ba_{\ba \cap \bb}X \mid X \in F^*_\ba \text{ is full over }\ba \cap \bb\}$$
Then $F^*$ is a proper $(\mu, A, \lambda)$-coherent filter.
\end{lemma}

{\bf Proof:} For each $\bb \in (\lambda)^{<\mu}$, $F^*_\bb$ is a proper filter by Lemma \ref{extfip-lem}.  Thus we must make sure the system is coherent.  Towards this, fix $\bb \subset \bc \in (\lambda)^{<\mu}$.

First, suppose that $X \in F_\bb^*$ and we want to show that $(\pi^\bc_\bb)^{-1} \in F_\bc^*$.  From the definition of $F^*_\bb$, there is $Y \in F_\bb$ and $X_1, \dots, X_n \in F^*_\ba$ full over $\ba \cap \bb$ such that 
$$Y \cap \bigcap_i (\pi^\bb_{\ba \cap \bb})^{-1} \pi^\ba_{\ba \cap \bb} X_i \subset X$$
Applying $(\pi^\bc_\bb)^{-1}$ to both sides and noting that $(\pi^{\bc}_{\bb})^{-1}(\pi^{\bb}_{\ba\cap\bb})^{-1}\pi^{\ba}_{\ba\cap\bb} = (\pi^{\bc}_{\ba\cap\bc})^{-1}\pi^{\ba}_
{\ba\cap\bc}$, we have 
$$(\pi^\bc_\bb)^{-1}Y \cap \bigcap (\pi^\bc_{\ba\cap \bc})^{-1} \pi^\ba_{\ba \cap \bc}X_i \subset (\pi^\bc_\bb)^{-1}X$$
Thus, $(\pi^\bc_\bb)^{-1}X \in F^*_\bc$, as desired.

Second, suppose that $X \in F_\bc^*$ is full over $\bb$ and we want to show $\pi^\bc_\bb X \in F^*_\bb$.  From the definition of $F^*_\bc$, there is $Y \in F_\bc$ and $X_1, \dots, X_n \in F^*_\ba$ full over $\ba \cap \bc$ such that 
$$Y \cap \bigcap_i (\pi^\bc_{\ba \cap \bc})^{-1} \pi^\ba_{\ba \cap \bc}X_i \subset X$$
Using Lemma \ref{nicefull-lem}, we can assume that $Y$ is full over $\bb$ and each $X_i$ is full over $\ba \cap \bb$; this implies $\pi^\ba_{\ba\cap \bc}X_i$ is full over $\ba \cap \bb$.  The goal is to show that 
$$\pi^\bc_\bb Y \cap \bigcap_i (\pi^\bb_{\ba \cap \bb})^{-1} \pi^\ba_{\ba \cap \bb} X_i \subset \pi^\bc_\bb X$$
Since $F$ is assumed to be a coherent filter, this suffices to show $\pi^\bc_\bb X \in F^*_\bb$.  Let $s \in \pi^\bc_\bb Y \cap \bigcap_i (\pi^\bb_{\ba \cap \bb})^{-1} \pi^\ba_{\ba \cap \bb} X_i$.  By definition, this means there are $t \in Y$ and $s_i \in X_i$ such that $s = \pi^\bc_\bb(t)$ and $\pi^\bb_{\ba \cap \bb} (s) = \pi^\ba_{\ba \cap \bb}(s_i)$.

We claim that $t \in (\pi^\bc_{\ba \cap \bc})^{-1} \pi^\ba_{\ba \cap \bc} X_i$.  To see this, we can compute that
\begin{eqnarray*}
\pi^{\ba\cap \bc}_{\ba \cap \bb}\left( \pi^\bc_{\ba \cap \bc}(t) \right) = \pi^{\bb}_{\ba\cap \bb}\left( \pi^\bc_\bb(t)\right) =\pi^\bb_{\ba \cap \bb}(s) = \pi^\ba_{\ba \cap \bb}(s_i) \in \pi^{\ba\cap \bc}_{\ba \cap \bb} \pi^\ba_{\ba \cap \bc} X_i
\end{eqnarray*}
Since $\pi^\ba_{\ba \cap \bc}X_i$ is full over $\ba \cap \bb$, this gives $\pi^\bc_{\ba\cap \bc}(t) \in \pi^\ba_{\ba \cap \bc}X_i$, as desired.

Thus, $t \in Y \cap \bigcap_i (\pi^\bc_{\ba\cap \bc})^{-1} \pi^\ba_{\ba \cap \bc}X_i$.  By assumption, we have $t \in X$ and $s = \pi^\bc_\bb(t) \in \pi^\bc_\bb X$, as desired.\hfill \dag\\

This brings us to the main theorem.

\begin{theorem}\label{cohultext-thm}
Any coherent filter can be extended to a coherent ultrafilter.
\end{theorem}

{\bf Proof:}  We use Zorn's Lemma.  Let $F$ be a coherent filter and $\mathcal{Z}$ be the collection of coherent filters extending $F$ ordered by $\Subset$.  Clearly, the union of any $\Subset$-increasing chain of coherent filters is a coherent filter, so by Zorn's Lemma there is a maximal element $E$.  If some $E_\ba$ is not an ultrafilter, then there is a proper filter $E^*_\ba$ extending it.  By Lemma \ref{onestep-lem}, this gives rise to $E^* \Supset E$, contradicting it's maximality.  Thus, $E$ is a coherent ultrafilter.\hfill \dag\\

Theorem \ref{cohultext-thm} suffices for Section 3, where we deal with definable coherent ultrapowers.  For completeness, we include some basic facts about coherent ultrafilters and coherent ultraproducts.  For the rest of the section, we focus on $(\mu, \kappa, \lambda)$-coherent ultrafilters.  Since we don't deal with definability, there is no loss.  First, note that Theorem \ref{cohultext-thm} does not prove the existence of proper coherent ultrafilter (although this follows by Theorem \ref{main-thm} applied to any proper extensions).  The following result says that we can generate a coherent ultrafilter starting with any collection of ``seeds" (at least when $\mu=\omega$).

\begin{theorem}
Let $\{U_\alpha \mid \alpha < \lambda\}$ be a collection of ultrafilters on $\kappa$.  Then there is an $(\omega, \kappa, \lambda)$-coherent ultrafilter $E$ such that, for all $\alpha < \lambda$, $E_{\{\alpha\}}$ is $U_\alpha$ after passing through the canonical bijection from ${}^1\kappa$ to $\kappa$.
\end{theorem}

For this proof, recall the product of filters: if $F_\ell$ is a filter on $I_\ell$ for $\ell = 0, 1$, then $F_0 \otimes F_1$ is the filter on $I_0 \times I_1$ given by
$$X \in F_0 \otimes F_1 \iff \left\{i \in I_0 \mid \{j \in I_1 \mid (i, j)\in X\} \in F_1 \right\} \in F_0$$
$F_0 \otimes F_1$ is a filter, is an ultrafilter iff $F_0$ and $F_1$ are, and the productive is associative, although non-commutative.  We use this product slightly modified to our situation so, e. g., the product of a filter on ${}^n\kappa$ and a filter on ${}^m\kappa$ is a filter on ${}^{n+m}\kappa$.

{\bf Proof:}  We define $E$ by induction by setting, for $\ba \in [\lambda]^{<\omega}$:
\begin{enumerate}
	\item If $\ba =\{\alpha\}\in[\lambda]^1$, then $E_{\ba} = U_{\alpha}$.
	\item If $\ba \in [\lambda]^n$ is $\alpha_0 < \dots < \alpha_{n-1}$, then $E_{\ba}$ is the common value of
	$$U_{\alpha_0} \otimes E_{\{\alpha_1, \dots, \alpha_{n-1}\}} = E_{\{\alpha_0, \dots, \alpha_{i-1}\}} \otimes U_{\alpha_i} \otimes E_{\{\alpha_{i+1}, \dots, \alpha_{n-1}\}} = E_{\{\alpha_0, \dots, \alpha_{n-2}\}} \otimes U_{\alpha_{n-1}}$$
	\end{enumerate}

$E$ is clearly a collection of ultrafilters, so it only remains to show coherence.  It is enough to check coherence for one-point extensions.

Let $\bb \in [\lambda]^{n+1}$ be $\alpha_0 < \dots < \alpha_n$, $X \subset {}^n\kappa$, and $i < n+1$.  Set $\ba = \bb -\{\alpha_i\}$ and   $X^* = \{s \in {}^{n+1} \kappa \mid \pi^{\bb}_{\ba}(s) \in X\}$.  For notational ease, set $\ba_1 =\{\alpha_0, \dots, \alpha_{i-1}\}$ and $\ba_2 = \{\alpha_{i+1}, \dots, \alpha_{n}\}$.  We have
\begin{eqnarray*}
X \in E_{\ba} &\iff& X \in E_{\ba_1} \otimes E_{\ba_2}\\
&\iff& \left\{s \in {}^i \kappa \mid \{ t \in {}^{n-i}\kappa \mid s^\frown t\in X\} \in E_{\ba_2}\right\} \in E_{\ba_1}\\
X^* \in E_{\bb} &\iff& X^* \in E_{\ba_1} \otimes U_{\alpha_i} \otimes E_{\ba_2}\\
&\iff& \left\{s \in {}^i \kappa \mid\left\{ j \in {}^1 \kappa\mid  \{ t \in {}^{n-i}\kappa \mid s^\frown j^\frown t\in X^*\} \in E_{\ba_2}\right\} \in U_{\alpha_i}\right\} \in E_{\ba_1}
\end{eqnarray*}
We have $\pi^{\bb}_{\ba}\left(s^\frown j^\frown t\right) = s^\frown t$, so $s^\frown t \in X$ iff $s^\frown j^\frown t \in X^*$.  This finishes the proof.\hfill\dag\\

Most applications of coherent ultrafilters (and extenders) involve taking the coherent relative of the ultrapower--which we call the coherent ultrapower--of a single model.  This is also true of our Theorem \ref{main-thm}.  However, one can naturally define a coherent ultraproduct.  In fact, there are at least two notions of a coherent ultraproduct one might define; one is more general than the other, but this extra generality borders on ``too general" and there seems to be no application of it that doesn't reduce to a simple case (yet).

First, the less general.  Let $E$ be an $(\mu, \kappa, \lambda)$-coherent ultrafilter, $\ba \in (\lambda)^\alpha$, and $\{M_s \mid s \in {}^\alpha\kappa\}$ be a collection of $\tau$-structures.  The \emph{coherent ultraproduct of $\{M_s \mid s \in {}^\alpha \kappa\}$ by $E$ at $\ba$} is denoted by $\prod^\ba M_s/E$ and is constructed as follows: for each $\bb\in {}^{\beta}\lambda$ that extends $\ba$, form the standard ultraproduct 
$$M^\bb_* := \prod_{s \in {}^\beta \kappa} M_{\pi^\bb_\ba(s)}/E_{\bb}$$
Then if $\bc \in (\lambda)^{<\mu}$ extends $\bb$ (and therefore also $\ba$), there is a natural map $f^{\bb, \bc}:M^{\bb}_* \to M^{\bc}_*$ by taking $[f]_{E_{\bb}}$ to $[f \circ \pi^{\bc}_{\bb}]_{E_{\bc}}$.  \L o\'{s}' Theorem shows that $f^{\bb, \bc}$ is an elementary embedding.  This is a directed system, and we take the colimit to form $\prod^\ba M_s/E$.  The coherent ultrapower of $M$ by $E$ is this construction with $M_s = M$.  This means that a coherent ultrapower is simultaneously a coherent ultraproduct at $\ba$ for every $\ba \in (\lambda)^{<\mu}$.  Note that the isomorphism type of the coherent ultraproduct is independent of the ordering on $\lambda$ (in the sense that any permutation of $\lambda$ induces an automorphism of coherent ultraproducts).

We can generalize this further by imposing the minimum structure necessary to make this construction work.  Fix an $(\mu, \kappa, \lambda)$-ultrafilter $E$ and structures $\{M^\ba_s \mid \ba \in (\lambda)^{<\mu}, s \in {}^{\otp(\ba)}\kappa\}$ such that if $\bb$ extends $\ba$ and $t \in {}^{\otp(\bb)}\kappa$, then $\tau\left(M^{\ba}_{\pi^\bb_{\ba}(t)}\right) \subset \tau(M^\bb_t)$ and there is elementary $f^{\ba, \bb}_{\pi^\bb_\ba(t), t}: M^\ba_{\pi^\bb_\ba(t)} \to M^\bb_t$.    Then we can again form the ultraproducts $M_*^\ba = \prod_{s \in {}^{|\ba|}\kappa} M^\ba_s/E_{\ba}$.  If $\bb$ extends $\ba$, then the $f^{\ba, \bb}_{s, t}$ induce $f^{\ba,\bb}_*:M^\ba_* \to M^\bb_*$ by taking $[f]_{E_{\ba}}$ to $[t \mapsto f^{\ba, \bb}_{\pi^\bb_\ba(s), s} \left( f(\pi^\bb_\ba(s))\right)]_{E_\bb}$.  Again, this is elementary by \L o\'{s}' Theorem and the elementarity of each $f^{\ba,\bb}_{s,t}$.  Set $\prod M^\ba_s/E$ to be the colimit of this system.

We can view coherent ultraproducts at some $\ba$ as these more general coherent ultraproducts by setting
$$M^\bb_s = \begin{cases} M_{\pi^\bb_\ba(s)} & \bb\text{ extends }\ba\\ \emptyset & \text{otherwise}\end{cases}$$
and each $f^{\ba, \bb}_{s,t}$ is the identity.  However, we know of no use for this extra level of generality.  As further evidence for their strangeness, suppose $E$ was a $(\kappa, \lambda)$-extender and $j_E:V \to M_E$ was the derived embedding.  Then an coherent ultraproduct at $\ba$ (for any $\ba \in [\lambda]^{<\omega}$) appears in $M_E$: if $f$ is the functions that takes $s \in {}^{|\ba|}\kappa$ to $M_s$, then $\prod M_s/E \cong j(f)\left(j^{-1}\rest j(\ba)\right)$.  However, the more general coherent ultraproducts don't seem to appear.  Some more work on the general properties of coherent ultraproducts can be found in \cite[Section 5]{b-mtlc}.

\section{The Definable Coherent Ultrapower}

We define the notion of a definable coherent ultrapowers, which combines the notion of a coherent ultrapower  and a definable ultrapower (see \cite{k-lw1w}).
\begin{defin}\label{frag-def}
    Fix a finitary language $\tau$ and a set $\mathcal{F} \subset \bL_{\infty,\omega}(\tau)$.
    \begin{enumerate}
        \item $\mathcal{F}$ is a \emph{fragment} iff it is closed under subformulae.
        \item $\mathcal{F}$ is an \emph{elementary fragment} iff it is a fragment closed under first-order operations, e.g., if $\phi(\bx,y) \in \mathcal{F}$, then $\exists y \phi(\bx,y)\in \mathcal{F}$.
    \end{enumerate}
\end{defin}

Note \cite[Chapter 4]{k-lw1w} and others use `fragment' to refer to what we call an `elementary fragment,' but it is a useful distinction.  It is easy to form fragments from sets of formulae.

\begin{prop}
Let $X \subset\bL_{\infty,\omega}(\tau)$.  Then there is a minimal fragment $\mathcal{F}_0$ containing $X$ and a minimal elementary fragment $\mathcal{F}_1$ containing $X$ with $\mathcal{F}_0 \subset \mathcal{F}_1$.
\end{prop}
We are going to work with the class $\K = (\Mod \psi, \prec_{\mathcal{F}})$.  Note that this includes first order classes (take $\mathcal{F}$ to be $\mathbb{L}_{\omega, \omega}(\tau)$ along with a single conjunction for the theory).  The following definitions are natural generalizations of those in \cite[Chapter 32]{k-lw1w}.

\begin{defin} \label{def-skol-def} Fix a sentence $\psi \in \mathcal{F}$.
\begin{enumerate}
	\item A $n$-ary \emph{$\psi$-definable function} $\phi(\bx, y)$ with domain $\chi(z)$ is a formula $\phi(\bx, y) \in \mathcal{F}$ such that 
	$$\psi \vDash \forall \bx \left( \bigwedge_{i<n} \chi(x_i) \to \exists^{=1} y\, \phi(\bx, y)\right)$$
	We sometimes write  $F_{\phi(\bx, y)}(\bx)$ for this function..
	\item Let $\chi(z) \in \mathcal{F}$.  We say that $\K$ has \emph{$\psi$-definable Skolem functions over $\chi(z)$} iff for every $\text{`}\exists y\, \phi(\bx, y)\text{'} \in \mathcal{F}$, there is a $\psi$-definable function $F^{\exists y \phi(\bx, y)}(\bx)$ with domain $\chi(z)$ such that 
	$$\psi \vDash \forall \bx \left( \bigwedge_i \chi(x_i) \wedge \exists y\, \phi(\bx, y) \to \phi\left( \bx, F^{\exists y' \phi(\bx', y')}(\bx)\right) \right)$$
	\item Suppose that $\K$ has $\psi$-definable Skolem functions over $\chi(z)$.  Let $M \in \K$ and $E$ be a $(\omega, \chi(M), \lambda)$-coherent ultrafilter for some $\lambda$.
	\begin{enumerate}
		\item Set $\prod^{\de}_{\chi,n} M$ to be set of all $n$-ary $\psi$-definable functions with parameters with domain $\chi(z)$.
		\item For $F, G \in \prod^{\de}_{\chi,n} M$ and $\ba \in [\lambda]^n$, set
		$$F \sim_{E_\ba} G \iff \{ s \in {}^n \chi(M) \mid M \vDash ``F(s) = G(s)"\} \in E_\ba$$
		\item Set $\prod^{\de}_{\chi,n} M /E_\ba$ to be the model with universe $\{ [F]_{E_\ba} : F \in \prod^{\de}_{\chi,n} M\}$ with the standard ultraproduct structure, e. g., if $R \in \tau$, then 
		$$R^{\prod^{\de}_{\chi,n} M /E_\ba}\left([F_1]_{E_\ba}, \dots, [F_n]_{E_\ba}\right) \iff \{ s \in {}^n \chi(M) \mid M\vDash R\left(F_1(s), \dots, F_n(s)\right) \} \in E_\ba$$
		\item Given $\ba \in [\lambda]^n$ and $\bb \in [\lambda]^m$ such that $\ba \subset \bb$, define $k_{\ba, \bb}: \prod^{\de}_{\chi, n} M/E_\ba \to \prod^{\de}_{\chi, m} M/E_{\bb}$ by
		$$k_{\ba, \bb}\left([F]_{E_{\ba}}\right) = \left[F \circ \pi^{\bb}_{\ba}\right]_{E_{\bb}}$$
		\item The definable coherent ultrapower of $M$ by $E$ on $\chi(z)$ is the directed colimit of the sequence
		$$\left\{\prod^{\de}_{\chi, n} M/E_\ba, k_{\ba, \bb} \mid \ba \subset \bb \in [\lambda]^{<\omega} \right\}$$
		And is denoted 
		$$\left( \prod^{\de}_{\chi} M/ E, k_{\ba, \infty}\right)_{\ba \in [\lambda]^{<\omega}}$$

	\end{enumerate}
	\item In the above definitions, if $\chi(z) \equiv \text{`}z = z\text{'}$, then we omit it, e. g., writing $\prod^\de M/E$.
\end{enumerate}
\end{defin}

We use $\mu = \omega$ here because there are no $\mathbb{L}_{\infty, \omega}$-definable functions with infinite arity; in turn, this is because formulae of $\bL_{\infty, \omega}$ definitionally have finitely many free variables (see \cite[Definition 4.9]{b-cofquant} for a use of a version of $\bL_{\infty, \omega}$ that allows infinitely many free variables).  However, this theory generalizes to the construction of definable $(\mu, \chi(M), \lambda)$-coherent ultraproducts in $\mathbb{L}_{\infty, \mu}$-axiomatizable classes (or $\mu$-AECs; see Boney, Grossberg, Lieberman, Rosicky, and Vasey \cite{bglrv-muaecs}, especially the Presentation Theorem 3.2 there).

We define the notion of Skolem functions for single quantifiers, but we note for later that we can also build Skolem functions for tuples.  Note that we crucially use that our fragment is closed under existential quantifiers in this result.

\begin{lemma}\label{nary-skol-lem}
Suppose that $\K$ has definable Skolem functions over $\chi$ and $\phi(\bx; y_1, \dots, y_\ell) \in \mathcal{F}$.  Then there are $\mathcal{F}$-definable functions $F_1(\bx), \dots, F_\ell(\bx)$ such that
$$\psi \vDash \forall \bx \left(\bigwedge_i\chi(x_i) \wedge \exists \by \phi(\bx, \by) \to  \phi\left(\bx, F_1(\bx), \dots, F_\ell(\bx)\right)\right)$$
\end{lemma}

{\bf Proof:}. We work by induction on $\ell = \ell(\by)$. $\ell=1$ is the definition of having definable Skolem functions.  Suppose this is true for some $\ell < \omega$ and we have $\phi(\bx; y_1, \dots,y_{\ell+1}) \in \mathcal{F}$.  By induction, we find $\mathcal{F}$-definable functions $F'_2(\bx, y_1), \dots, F'_{\ell+1}(\bx, y_1)$ for those variables.  Then we apply the definition of definable Skolem functions for 
$$\exists y_1\phi\left(\bx, y_1, F'_2(\bx, y_1), \dots,F'_{\ell+1}(\bx, y_1)\right)$$
to get the function $F_1(\bx)$.  Then, for $k>1$, we define
$$F_k(\bx):=F'_k(\bx, F_1(\bx))$$
\hfill\dag\\
During Definition \ref{def-skol-def}, we assumed that certain definitions were well-defined.  We note this now.

\begin{prop}
The construction of $\prod^\de_\chi M/E$ is well defined; that is,
\begin{enumerate}
	\item $\sim_{E_\ba}$ is a $\tau$-congruence relation on $\prod^\de_{\chi, n} M$; and
	\item given $\ba \subset \bb \subset \bc \in [\lambda]^{<\omega}$, we have that $k_{\ba, \bb}$ is a $\tau$-embedding $k_{\ba, \bc} = k_{\bb, \bc} \circ k_{\ba, \bb}$.
\end{enumerate}
\end{prop}

{\bf Proof:} The proof is a straightforward calculation.\hfill \dag\\

It will be useful to recall the form of any element in a $(\mu, \chi(M), \lambda)$-coherent ultrapower: given $x \in \prod^\de_\chi M/E$, we have $x = k_{\ba, \infty}\left([F]_{E_\ba}\right)$ for some $\ba$ and $F$, where $F$ is a $|\ba|$-ary definable function from $\chi(M)$ to $M$; in fact, there are many such $\ba$ and $F$.  We will write this as $[\ba, F]_E$.  Then it is easy to check that
$$[\ba,F]_E = [\bb, G]_E \iff \left\{ s \in {}^{|\ba \cup \bb|}\chi(M) \mid M \vDash \text{`}F \circ \pi^{\ba\cup \bb}_\ba(s) = G \circ \pi^{\ba \cup \bb}_\bb(s)\text{'} \right\} \in E_{\ba \cup \bb}$$

One final notion is necessary to show that this construction interacts well with the (potentially) non-elementary class $\K$.  Typically, highly complete (and therefore non-ZFC) ultrafilters are necessary to preserve $\mathbb{L}_{\infty, \omega}$ formulas.  However, since we only deal with \emph{definable} functions, we have more leeway.

\begin{defin} \
\begin{enumerate}
	\item Suppose $M \in \K$, $\phi(x_1, \dots, x_n) \in \mathcal{F}$, and $G_1, \dots, G_m$ are $n$-ary $\psi$-definable functions with domain $\chi(M)$.  Then set
	$$\phi \circ (G_1, \dots, G_m) \left(\chi(M)\right) := \{s \in {}^n \chi(M) \mid M \vDash\text{`} \phi\left(G_1(s), \dots, G_m(s)\right)\text{'} \}$$
	\item For $M \in \K$, we say that a filter $F$ on ${}^n \chi(M)$ is \emph{$\mathcal{F}$-complete} iff for all $\text{`}\bigwedge_{\alpha < \kappa} \phi_\alpha(x_1, \dots, x_m)\text{'} \in \mathcal{F}$ and $n$-ary $\psi$-definable functions $G_1, \dots, G_m$ with domain $\chi(M)$, we have
    $$\text{for all $\alpha < \kappa$}, \phi_\alpha \circ(G_1, \dots, G_m) \left(\chi(M)\right) \in F \implies \bigwedge_{\alpha < \kappa}\phi_\alpha \circ(G_1, \dots, G_m) \left(\chi(M)\right) \in F$$
	\item We say that a coherent filter $F$ is $\mathcal{F}$-complete iff each $F_\ba$ is.
\end{enumerate}
\end{defin}

\begin{prop} \label{basicprop-prop}
Suppose $\K$ has definable Skolem functions over $\chi(z)$, $M \in \K$, and $E$ is an $\mathcal{F}$-complete $(\omega, \chi(M), \lambda)$-coherent ultrafilter.
\begin{enumerate}
    \item\label{elemfrag} If $\mathcal{F}^*$ is the minimal elementary fragment containing $\mathcal{F}$, then $E$ is $\mathcal{F}^*$-complete.
	\item \label{elememb} For all $\ba \subset \bb \in[\lambda]^{<\omega}$, the embeddings $k_{\ba, \bb}$ and $k_{\ba, \infty}$ are $\mathcal{F}$-elementary.
	\item \label{cohlos} Let $[\ba_1, G_1]_E, \dots, [\ba_m, G_m]_E\in \prod^\de_\chi M/E$ and $\phi(x_1, \dots, x_m) \in \mathcal{F}$.  Then
	$$\prod^\de_\chi M/E \vDash \phi\left([\ba_1, G_1]_E, \dots, [\ba_n, G_m]_E\right)$$
	$$\text{iff}$$
	$$\phi \circ(G_1 \circ \pi^{\cup \ba_i}_{\ba_1}, \dots, G_m \circ \pi^{\cup \ba_i}_{a_n}) \left(\chi(M)\right) \in E_{\cup \ba_i}$$
	\item \label{closed} $\prod^\de_\chi M/E \in \K$
	\item \label{power} The coherent ultrapower embedding $j_E: M \to \prod^\de_\chi M/E$ is $\mathcal{F}$-elementary and is also proper iff some $E_\ba$ is non-principal.
\end{enumerate}
\end{prop}

We will refer to item (\ref{cohlos}) as \L o\'{s}' Theorem for definable coherent ultraproducts.  Note that the set appearing there is 
$$\left\{s \in{}^{|\cup\ba_i|} \chi(M) \mid M \vDash \phi\left( G_1(\pi^{\cup \ba_i}_{\ba_1}(s)), \dots, G_n(\pi^{\cup \ba_i}_{\ba_n}(s)) \right) \right\}$$
In the proof, we will use the fact that, if $G$ is a definable function, then so is $G \circ \pi^\bb_\ba$.

{\bf Proof:}  \ref{elemfrag}) Straight from the definition since $\mathcal{F}$ and $\mathcal{F}^*$ contain the same formulae of the form
$$\text{`}\bigwedge_{\alpha < \kappa} \phi_\alpha(x_1, \dots, x_m)\text{'}$$\\

\ref{elememb}) We will use \L o\'{s}' Theorem for definable ultraproducts, see \cite[Theorem 46]{k-lw1w}.

Suppose that $\ba \subset \bb \in [\lambda]^{<\omega}$ and $[G_1]_{E_{\ba}}, \dots, [G_m]_{E_{\ba}} \in \prod^\de_{\chi, n} M/E_\ba$.  Let $\phi(x_1, \dots, x_m) \in \mathcal{F}$.  Then
\begin{eqnarray*}
\prod^\de_{\chi, n} M/E_\ba \vDash \phi\left([G_1]_{E_{\ba}}, \dots, [G_m]_{E_{\ba}}\right) &\iff& \phi\circ(G_1, \dots, G_m)\left(\chi(M)\right) \in E_\ba\\
&\iff& (\pi^\bb_\ba)^{-1}\phi\circ(G_1, \dots, G_m)\left(\chi(M)\right)\\
&& =\phi\circ(G_1 \circ \pi^\bb_\ba, \dots, G_m\circ \pi^\bb_\ba)\left(\chi(M)\right) \in E_\bb\\
&\iff& \prod^\de_{\chi, m} M/E_\bb \vDash \phi\left(k_{\ba, \bb}([G_1]_{E_{\ba}}), \dots, k_{\ba, \bb}([G_m]_{E_{\ba}})\right)
\end{eqnarray*}
We used the coherence of $E$ in the second line, and the definition of $k_{\ba, \bb}$ in the final line.  The $\mathcal{F}$-elementarity of the $k_{\ba, \infty}$ follows from the $\mathcal{F}$-elementarity of the $k_{\ba, \bb}$ since $\mathcal{F}$-elementary maps are closed under directed colimits\footnote{Importantly, this does not hold for $\bL_{\omega_1, \omega_1}$-elementary maps, which is why we have restricted to fragments of $\bL_{\infty, \omega}$. A simple example to demonstrate this is the property of being well-founded: the ill-founded order $(\mathbb{Z}, <)$ is the union of well-founded $\left([-n,n], <\right)$}.\\

\ref{cohlos})  This follows from (\ref{elememb}) and the fact that $[\ba_i, G_i]_E = k_{\ba, \infty}([G_i]_{E_{\ba_i}}) = k_{\cup \ba_j, \infty}([G_i \circ \pi^{\cup \ba_j}_{\ba_i}]_{E_{\cup \ba_j}}$.  Then 
\begin{eqnarray*}
& &\prod^\de_\chi M/E \vDash \phi\left([\ba_1, G_1]_E, \dots, [\ba_n, G_m]_E\right) \\
&\iff& \prod^\de_{\chi, n} M /E_{\cup\ba_i} \vDash \phi \left([G_1 \circ \pi^{\cup\ba_i}_{a_1}]_{E_{\cup a_i}}, \dots, [G_m \circ \pi^{\cup \ba_i}_{\ba_n}]_{E_{\cup\ba_i}} \right) \\
&\iff& \phi \circ (G_1 \circ \pi^{\cup \ba_i}_{\ba_1}, \dots, G_m \circ \pi^{\cup \ba_i}_{\ba_n}) \left(\chi(M)\right) \in E_{\cup \ba_i}
\end{eqnarray*}

\ref{closed}) $\psi \in \mathcal{F}$, so this follows from (\ref{elememb}).\\

\ref{power}) For $\ba \in [\lambda]^{<\omega}$ and $p \in M$, set $G^{\ba}_p$ to be the constant $p$-valued function with domain ${}^\ba \chi(M)$; this is definable with the parameter $p$.  Then $[\ba, G^\ba_p]_E = [\bb, G^\bb_p]_E$ for all $\ba, \bb \in [\lambda]^{<\omega}$.

Define $j_E$ to take $p \in M$ to $[\ba, G^\ba_p]_E$ for some/any $\ba \in [\lambda]^{<\omega}$.  Then $j_E$ is elementary by \L o\'{s}' Theorem.

For the properness, first suppose that $[\ba, H]_E \in \prod^\de_\chi M/E$ is not in the range of $j_E$.  If $E_\ba$ is principal, then it is generated by some $s \in {}^{|\ba|} \chi(M)$.  Then $[\ba, H]_E = [\ba, G^\ba_{H(s)}]_E$, contradicting $[\ba, H]_E$ not being in the range of $j_E$.  Second, suppose that some $E_\ba$ is non-prinicipal.  Let $\id_\ba$ be the identity function on ${}^{|\ba|}\chi(M)$.  Then $[\ba, \id_\ba]_E \not\in j_E(M)$.\hfill \dag\\

This shows that the existence of an $\mathcal{F}$-complete coherent ultrafilter will give rise to an extension.  Our main theorem is that the converse holds as well.

\begin{theorem} \label{main-thm}
Suppose that $\K$ has definable Skolem functions over $\chi(z)$ and $M \prec_\mathcal{F} N$ where $M, N \in \K$.  Let $\lambda = |\chi(N)|$ and $f:\lambda\to \chi(N)$ be a bijection. Then there is an $\mathcal{F}$-complete $(\omega, \chi(M), \lambda)$-coherent ultrafilter $E$ and $\mathcal{F}$-elementary $h: \prod^\de_\chi M/E \to N$ such that $h \circ j_E = \id_M$ and 
$$h"\chi\left(\prod^\de_\chi M/E\right) = \chi(N)$$
\end{theorem}

{\bf Proof:}  By Proposition \ref{basicprop-prop}.(\ref{elemfrag}), we may assume that $\mathcal{F}$ is elementary.  For $\ba \in [\lambda]^n$, define $F_\ba$ to be the filter generated by sets of the form $\phi\circ(G_1, \dots, G_m)\left(\chi(M)\right)$ for $\phi(x_1, \dots, x_m) \in \mathcal{F}$ and $G_1, \dots, G_m \in\prod^\de_{\chi, n} M$ such that 
$$N \vDash \phi\left(G_1(f"\ba), \dots, G_m(f"\ba)\right)$$
Since $\mathcal{F}$ is closed under finitary conjunctions,  $F_{\ba}$ consists of all $X \subset {}^{\otp(\ba)}\chi(M)$ that contain one of these sets.
\begin{claim}\label{claim1}
$F$ is an $\mathcal{F}$-complete $(\omega, \chi(M), \lambda)$-coherent filter.
\end{claim}

{\bf Proof:} Each $F_\ba$ is a filter by definition (and will be principal exactly when $f"\ba \in \chi(M) \subset \chi(N)$).  Nearly by definition, each $F_\ba$ is $\mathcal{F}$-complete.

To show coherence, let $\ba \subset \bb \in [\lambda]^{<\omega}$.  We use the characterization following Definition \ref{main-def}.

First, if $X \in F_\ba$, then it contains
$$\phi\circ\left(G_1, \dots, G_{m}\right)\left(\chi(M)\right)$$
such that
$$N \vDash \phi\left(G_1(f"\ba), \dots, G_{m}(f"\ba)\right)$$
Then $(\pi^\bb_\ba)^{-1}X$ contains 
$$\phi\circ (G_1 \circ \pi^\bb_\ba, \dots, G_{m} \circ \pi^\bb_\ba)\left(\chi(M)\right)$$
and, because $f"\pi^{\bb}_{\ba}(\bb) = f"\ba$,
$$N \vDash \phi \left(G_1\left(f"\pi^\bb_\ba(\bb)\right), \dots, G_{m}\left(f"\pi^\bb_\ba(\bb)\right)\right)$$
Thus, $(\pi^\bb_\ba)^{-1}X \in F_\bb$.

Second, suppose $X \in F_\bb$ is full over $\ba$ and it contains $\phi \circ(G_1, \dots, G_{m})\left(\chi(M)\right)$ with
$$N \vDash \phi(G_1(f"\bb), \dots, G_m(f"\bb))$$
We can write each $G_j$ as $G_j(\bx; \by)$, where $\bx$ corresponds to $\ba$ and $\by$ corresponds to $\bb-\ba$.  We can use Lemma \ref{nary-skol-lem} to find definable Skolem functions $F_1(\bx), \dots, F_k(\bx)$ where $k = \ell(\by)$.

Then set 
$$G^{*}_j(\bx)=G_j\left((\bx,F_1(\bx), \dots, F_\ell(\bx)\right)$$
These definable functions have the property that
$$\psi\vDash\forall\bx \left(\bigwedge\chi(x_j) \wedge\exists\by \phi(G_1(\bx,\by), \dots, G_{m}(\bx, \by) \to \phi\left((G^{*}_1(\bx), \dots, G^{*}_{m}(\bx)\right)\right) $$
Note that $G_j^{*} \in \prod^\de_{\chi, \ell(\bx)} M$.

\begin{claim}\label{claim2}
$\phi \circ (G_1^{*}, \dots, G_{m}^{*}) \left(\chi(M)\right) \subset \pi^\bb_\ba X$
\end{claim}

{\bf Proof:} Let $s \in\phi \circ (G_1^{*}, \dots, G_{m}^{*}) \left(\chi(M)\right)$ and set $t = F_1^M(s), \dots, F_\ell^M(s)$.  Thus, $M \vDash \text{`}\phi\left(G_1(s;t), \dots, G_{m}(s; t)\right)\text{'}$.  So
$$s{}^\frown t \in \phi \circ(G_1, \dots, G_{m})\left(\chi(M)\right) \subset X$$
$$s = \pi^\bb_\ba(s{}^\frown t) \in \pi^\bb_\ba X$$
 \hfill $\dag_{\text{Claim \ref{claim2}}}$\\
By assumption that $X\in F_\bb$ and the condition on definable Skolem functions, we have
\begin{eqnarray*}
& & N \vDash \phi(G_1(f"\ba; f"(\bb-\ba)), \dots, G_m(f"\ba;f"(\bb-\ba)))\\
&\implies& N \vDash \phi(G_1^*(f"\ba), \dots, G^*_m(\ba))
\end{eqnarray*}
Then $\pi^\bb_\ba X \in F_\ba$. \hfill $\dag_{\text{Claim \ref{claim1}}}$\\

Now, by Theorem \ref{cohultext-thm}, there is an $(\omega, \chi(M), \lambda)$-coherent ultrafilter $E$ extending $F$.  Since all of the relevant sets--the $\phi \circ (G_1, \dots, G_n)\left(\chi(M)\right)$--are already in $F$, the $\mathcal{F}$-completeness of $F$ transfers to $E$.

Thus, we can form $\prod^\de_\chi M/E$ and $\mathcal{F}$-elementary $j_E:M \to \prod^\de_\chi M/E$ by Proposition \ref{basicprop-prop}.  We build a system of maps $\{h_\ba:\prod^\de_{\chi, n} M/E_\ba \to N \mid \ba \in [\lambda]^{<\omega}\}$ as in \cite[Corollary following Theorem 47]{k-lw1w}: set $h_\ba\left([G]_{E_\ba}\right) = G^N(f"\ba)$.  The same argument as \cite{k-lw1w} shows that this is a well defined $\mathcal{F}$-embedding.  It is also a straightforward calculation to show that, for any $\ba \subset \bb \in [\lambda]^{<\omega}$, 
$$h_\ba = h_\bb \circ k_{\ba, \bb}$$
Then we can push the maps through the colimit to get $\mathcal{F}$-elementary $h_{\infty}: \prod^\de_\chi M/E \to N$ such that $h_\ba = h_\infty \circ k_{\ba, \infty}$; this has an explicit description
$$h_\infty([\ba, G]_E) = G^N(f"\ba)$$

Finally, we show that $h_\infty$ is surjective.  Let $n \in \chi(N)$.  Clearly, the identity function is definable, so
$$h_\infty\left([\{f^{-1}(n)\}, \id]_E\right) = n$$
 \hfill $\dag_{\text{Theorem \ref{main-thm}}}$\\

This gives a nice corollary that characterizes $\mathcal{F}$-elementary extensions by coherent ultrafilters.

\begin{cor}\label{elmextchar-cor}
Suppose that $\K$ has definable Skolem functions and $M \prec_{\mathcal{F}} N$ from $\K$.  Then there is an $\mathcal{F}$-complete $(\omega, |M|, \|N\|)$-coherent ultrafilter $E$ such that $N$ and $\prod^\de M/E$ are isomorphic via a map extending $j$.
\end{cor}

Note that this is stronger than related results in set theory, e. g., an extender derived from a strong cardinal can only capture an initial segment of the target model, while Corollary \ref{elmextchar-cor} captures the entire target model.  This comes from the fact that our target model $N$ is a set, while the target model of a strong embedding is a proper class.

\section{Examples} \label{ex-sec}

\subsection{First-order}

We apply this theory to some elementary classes with definable Skolem functions.  Here, $\mathcal{F}$ can be taken to be $\mathbb{L}_{\omega, \omega}$ (or, pedantically, the smallest fragment of $\mathbb{L}_{|T|^+, \omega}$ containing $\bigwedge_{\phi \in T} \phi$, but those are equivalent).  This simplifies matters as the property of a coherent ultrafilter being $\mathcal{F}$-complete is vacuously satisfied.

Beyond the ones listed below, other examples of first order theories with definable Skolem functions include $p$-adically closed fields and real closed domains \cite[Theorem 3.2 and Corollary 3.6]{d-defskol}.  A more complete look at which $o$-minimal theories have definable Skolem functions can be found in \cite{de-defskol}.

\subsubsection{Peano Arithmetic}
$PA$ has definable Skolem functions.  $\mathbb{N}$ is a prime model of $PA$ \cite[Example 3.4.5]{changkeisler}, but in the sense of $\subset$ rather than $\prec$.  If we move to True Arithmetic ($TA=Th(\mathbb{N})$), then $\mathbb{N}$ is prime in the desired sense.  This means that every model of $TA$ can be viewed (up to isomorphism) as an extension of $\mathbb{N}$ and, using Theorem \ref{main-thm}, can be characterized by a $(\omega, \omega, \lambda)$-coherent ultrafilter.

Let $M \succ \mathbb{N}$.  Fix some bijection $f:\lambda \to |M|$.  Define $F^M$ to be the $(\omega, \omega, \lambda)$-coherent filter defined by, for $\ba \in [\lambda]^n$, $F^M_\ba$ is all subsets of $\mathbb{N}$ that contains some $\phi(\mathbb{N})$ such that $M \vDash \phi(f"\ba)$.  This can be extended to a $(\omega, \omega, \lambda)$-coherent ultrafilter $E$.  Define $h:M \to \prod^\de \mathbb{N}/E$ as in Theorem \ref{main-thm} to take $m \in M$ to $[\{f^{-1}(m)\}, \id]_E$.  This is an isomorphism.

In general, any consistent completion of $PA$ has some prime model $M$ and a similar analysis can be done for those models.

\subsubsection{$o$-minimal Theories and Transseries}

A more complex example can be made from any $o$-minimal theory that expands a group since all such theories have definable Skolem functions \cite{d-ominbook}.  As an application of this, one could take the restriction of transseries $\mathbb{T}$ (see \cite{avv-trans}) to the language of ordered fields and view this as a extension of $\mathbb{R}$.  Then, we can write $\mathbb{T}$ as a definable ultrapower of $\mathbb{R}$, which encodes it into the combinatorial of a coherent ultrafilter $E^{\mathbb{T}}$ on $\mathbb{R}$. \\

\subsubsection{Extensions of $ZFC$} In models of set theory, definable Skolem functions are intimately connected with the inner model $HOD$ of hereditarily ordinal definable sets.  $V=HOD$ is equivalent to the assertion that there is a definable well-ordering of $V$, which immediately gives definable Skolem functions.  In fact, a completion of $ZF$ has definable Skolem functions iff it proves `$\exists x \left(V = HOD(x)\right)$.'  In these cases, we can use Theorem \ref{main-thm} to understand the elementary extensions of these models as above.  Contrasting this, there are models of set theory beyond parameterized $HOD$ that still support Skolem functions over some set; for an example, see the notion of `$M$ satisfies choice over $x$' in \cite[Chapter 33]{k-lw1w}.

\subsection{Non-example: Skolemizing first-order theories} One might hope to apply Theorem \ref{main-thm} to any first-order theory by first Skolemizing it\footnote{The 'one' in this sentence is a reference to a past version of the author who was unaware of Keisler's problem and the difficulties outlined in this discussion.  Thanks again to the referee for pointing this out!}.  In order to do so, one would need to take a pair $M \prec N$ and lift it to Skolemized expansions $M^{sk} \prec N^{sk}$.  However, this expansion might not always exist.  Indeed, in 1967-68, Keisler asked if it was always possible to do even for a single Skolem function.  Over the years, Enayat, Keisler, Knight, Lachlan, Payne, and Winkler \cite{e-tree-keisler, k-expansions, k-skolem, l-skolem, p-skolem, w-skolem} have given various negative and positive examples to Keisler's problem (\cite{e-tree-keisler} gives a longer history).  If we are able to expand, we can show the following.

\begin{cor}
Let $T$ be a first order theory, $M \prec N$ models of $T$ with Skolem expansions so $M_{Sk} \prec N_{Sk}$. Then there is an $(\omega, |M|, \|N\|)$-coherent ultrafilter $E$ such that $\left(\prod^\de M_{Sk}/E \right)\rest \tau(M)$ and $N$ are isomorphic over $M$.
\end{cor}

However, the Skolemization requires extra information and $\prod^\de M_{Sk}/E \rest \tau(M)$ is very different from $\prod^\de M/E$ (the latter of which might not be a structure).  Thus we examine a few cases with definable Skolem functions.

\subsection{Beyond first-order}

Now we explore a general method of applying this construction to nonelementary classes.  We use Abstract Elementary Classes (AECs) as our framework of choice since they encompass many others.  This is not always ideal (the `Skolem functions' from Shelah's Presentation Theorem are not always the best possible), and we discuss other cases below.  For this section, we assume that the reader is familiar with AECs and point them towards Baldwin \cite{baldwinbook} as a reference.

Let $\K = (\K, \prec_\K)$ be an AEC with $\kappa = \LS(\K)$.  Given $M \prec_\K N$, we would like to realize $N$ as a definable coherent ultrapower of $M$.  However, we have no notion of definability, so this is difficult (to say the least).  Thus, we must appeal to Shelah's Presentation Theorem \cite[Conclusion 1.13]{sh88} (see also \cite[Theorem 4.15]{baldwinbook}).  This says that we can represent $\K$ as follows:
\begin{itemize}
	\item $\tau^*:= \tau(\K) \cup \{ F_n^i(\bx) \mid n < \omega, i < \kappa\}$;
	\item $T^* = \{ \forall \bx F_n^k(\bx) = x_k \mid k < n < \omega\}$ is a (first-order) $\tau^*$-theory;
	\item $\Gamma$ is a set of quantifier-free $\tau^*$-types;
	\item the models of $\K$ are precisely the $\tau(\K)$-reducts of models of $T^*$ omitting each type in $\Gamma$ (this is written $PC(T^*, \Gamma, \tau(\K))$); and
	\item for $M, N \in \K$, we have $M \prec_\K N$ iff there are expansions $M^*$ and $N^*$ of $M$ and $N$ to models of $T^*$ that omit $\Gamma$ and $M^* \subset N^*$.
\end{itemize}
Set 
$$\psi : = \bigwedge_{p \in \Gamma} \forall x_1, \dots, x_n \bigvee_{\phi \in p} \neg \phi(\bx)$$
and $\mathcal{F}$ to be the smallest fragment of $\mathbb{L}_{(2^\kappa)^+, \omega}(\tau^*)$ containing $\psi$ that is closed under finite conjunctions.  Note that the definable functions of this language are just the functions of the language and that it trivially has definable Skolem functions (because $\mathcal{F}$ has no existentials).

Now suppose that $M \prec_\K N$.  Then we can expand these to $\tau^*$-structures satisfying $\psi$ such that $M^* \subset N^*$.  From the definition of the fragment, this gives $M^* \prec_{\mathcal{F}} N^*$.  Thus, we can use Theorem \ref{main-thm} to find an $\mathcal{F}$-complete $(\omega, |M|, \|N\|)$-coherent ultrafilter $E$ such that $\prod^{def} M^* /E \cong N^*$.  Restricting to the original language, we have 
$$\left(\prod^{def} M^* /E\right) \rest \tau(K) \cong N$$

In general, we have that the ultraproduct operation commutes with the restriction of languages.  However, this is not the case with \emph{definable} ultraproducts because the notion of whether or not a function is definable is very language dependent.

The expansion to the Shelah Presentation is very coarse in that it pays no attention to how the original AEC was defined.  We might have even started with a very nice first-order class with built-in Skolem functions, and then made it more complex by adding extra functions.  Working with the original definition can often lead to nicer Skolem functions.\\

Recall that $\mathbb{L}(Q)$ is first-order logic with an additional quantifier $Q$ that stands for ``there exists uncountably many.''  If $T$ is a $\mathbb{L}(Q)(\tau)$-theory and $\mathcal{F} \subset \mathbb{L}(Q)(\tau)$ is a fragment containing $T$, then we  define $\prec^*_\mathcal{F}$ by
\begin{eqnarray*}
M \prec_\mathcal{F}^* N &\iff& M \prec_{\mathcal{F}} N \text{ and if }\ba \in M \text{ and }M \vDash Qx\phi(x, \ba)\text{,}\\
& &\text{ then } \phi(N, \ba) = \phi(M, \ba)
\end{eqnarray*}
where $\prec_{\mathbb{L}(Q)}$ is elementary substructure according to the logic, then $\K := \left(\Mod T, \prec^*\right)$ is an AEC with $\LS(\K) = \aleph_1 + |T|$ (see \cite[Exercise 5.1.3]{baldwinbook}).  We wish to add Skolem function to this class.  One option is to use Shelah's Presentation Theorem.  A more precise option is to do the following:
\begin{itemize}
	\item $\tau_{Sk} = \tau \cup \{F_{\phi}(\by) \mid \exists x \phi(x, \by) \in \mathcal{F}\} \cup\{G^{0,\alpha}_{\phi}(\by) \mid Qx \phi(x, \by) \in \mathcal{F}, \alpha < \omega_1\} \cup\{ G^{1,n}_{\phi}(\by) \mid \neg Qx \phi(x, \by) \in \mathcal{F}, n < \omega\}$
	\item $T_{Sk}$ is $T$ with the following changes:
	\begin{enumerate}
		\item every place ``$\exists x \phi(x, \by)$'' appears, replace it with ``$\phi(F_\phi(\by), \by)$''
		\item every place ``$Qx \phi(x, \by)$'' appears, replace it with ``$\wedge_{\alpha < \omega_1} \phi(G^{0, \alpha}_{\phi}(\by), \by)$"
		\item every place ``$\neg Q x \phi(x, \by)$" appears, replace it with ``$\forall z (\phi(z, \by) \to \vee_{n < \omega} z = G^{1,n}_\phi(\by)$"
	\end{enumerate}
	\item Let $\mathcal{F}_{Sk}$ be the smallest fragment of $\mathbb{L}_{\omega_2, \omega}(\tau)$ containing $T_{Sk}$.
\end{itemize}

Then $\K_{Sk} := (\Mod T_{Sk}, \prec_{\mathcal{F}_{Sk}})$ is the Skolemization of $\K$ and the above results can be applied.  Note that this Skolemization is $\mathbb{L}_{\infty, \omega}$-axiomatizable even though $\mathbb{L}(Q)$ is only $\mathbb{L}_{\omega_1, \omega_1}$-axiomatizable.

Similar refined Skolemizations for AECs defined using `existentially defined' quantifiers are found in \cite[Section 4]{b-cofquant}.

\bibliographystyle{amsalpha}
\bibliography{bib}

\end{document}